\documentclass[12pt]{amsart}

\usepackage{amssymb}
\newtheorem{prop}{Proposition}
\newtheorem{lemma}[prop]{Lemma}

\newtheorem{thm}[prop]{Theorem}

\theoremstyle{definition}

\newcommand{\Fp}{\mathbb{F}_p}
\begin{document}
\title[Residues in approximate
subgroups]{Distribution of residues in approximate subgroups of $\Fp^*$}


\author[N. Hegyv\'ari]{Norbert Hegyv\'ari}
\address{Norbert Hegyv\'{a}ri, ELTE TTK,
E\"otv\"os University, Institute of Mathematics, H-1117
P\'{a}zm\'{a}ny st. 1/c, Budapest, Hungary}
\email{hegyvari@elte.hu}

\author[F. Hennecart]{Fran\c cois Hennecart}
\address{Fran\c cois Hennecart,
PRES Universit\'e de Lyon,
Universit\'e Jean-Monnet,
Laboratoire de math\'ematiques de l'Universit\'e de Saint-\'Etienne
23, rue Michelon,
42023 Saint-\'Etienne, France} \email{francois.hennecart@univ-st-etienne.fr}

\thanks{Research of the first author is partially supported by OTKA grants K~61908, K 67676}
\thanks{The first author is grateful to the members of the
LAMUSE (Laboratory of Mathematics of the  University of Saint-Etienne) for their warm hospitality during his stay}

\date{\today}

\begin{abstract}
We extend a result due to Bourgain on the uniform distribution of residues by proving that subsets
of the type $f(I)\cdot H$ is equidistributed
(as $p$ tends to infinity) where $f$ is a polynomial,
$I$ is an interval of $\Fp$  and
$H$ is an approximate subgroup of $\Fp^*$ with size larger  than polylogarithmic in $p$.
\end{abstract}

\maketitle

\section{\bf Introduction}

Since few decades, additive combinatorics has become a central topic in number theory. At the origin,  there are several
very powerful and important results such as Freiman's theorem, Szemer\'edi's theorem,
Balog-Szemer\'edi-Gowers' theorem, etc.
(see \cite{TV} for a detailed  description of these results). Surprisingly, these results have many applications not only in combinatorial additive number theory but also in various topics such as the estimation of exponential sums. In this paper, we consider the closely related question of the equidistribution of the elements of a given multiplicative subgroup of a finite field
with prime cardinality. For $\delta$ a positive real number and $g$ a non zero  element of the prime field with $p$ elements, Bourgain obtained
in \cite{B} , under some restricted condition on the order of $g$,
an asymptotic equidistribution for the residues  $xg^n\pmod{p}$, $0\le x<p^{\delta},\ n\ge0$.
The proof of Bourgain's result uses the three above quoted theorems in combination with algebraic tools.
In this paper, our aim is to extend this equidistribution result to less structured  sets.

For any prime number $p$, we denote by $\Fp$
the finite field with $p$ elements, and let
$\Fp^*=\Fp\smallsetminus\{0\}$.
By an interval in $\Fp$, we mean a subset in the form
$$
I=\{ax+b\pmod{p}\;\mid\, 0\le x\le |I|-1\}
$$
for some $(a,b)\in \Fp^{*}\times \Fp$.

Let $\epsilon>0$. We say that a subset $H$ of $\Fp$ is
$\epsilon$-equidistributed modulo $p$ if for any
interval $J$ of $\Fp$, we have
$$
\Big|\frac{|H\cap J|}{|H|}-\frac{|J|}{p}\Big|<\epsilon.
$$
Using the Weyl criterion, $\epsilon$-equidistribution modulo $p$
follows from the following bound on trigonometric sums
\begin{equation}\label{eqn1}
\max_{\substack{a\in\mathbb{Z}\\p\,\nmid\, a}}\Big|\sum_{h\in H}
e_p(ah)\Big|<\epsilon' |H|,
\end{equation}
for some  $\epsilon'>0$ (depending on $\epsilon$).
Here we wrote $e_p(x)$ for $e^{2\pi i x/p}$.

Let $H$ be a (multiplicative) subgroup of $\Fp^*$.
In \cite{BGK}, it is shown that if
$$
|H|>p^{(\log\log p)^{-c}}
$$
where $c$ is an explicit positive constant then
\eqref{eqn1} holds true if $p$ is large enough.

In another direction, Bourgain considered in \cite{B0}
the question of the distribution of the values taken by
a \textit{sparse} polynomial formed by monomials
wich are, in some sense, sufficiently independent:
for $a_1,\dots,a_d\in\Fp^*$,
$$
\gcd(k_i,p-1),\gcd(k_i-k_j,p-1)<p^{1-\gamma}\Longrightarrow\Big|\sum_{x=1}^{p}e_p(a_1x^{k_1}+\cdots+a_dx^{k_d})\Big|
<p^{1-\delta},
$$
where $\delta=\delta(\gamma)>0$.

Bourgain also investigated in \cite{B} a mixed question and showed that $I\cdot H:=\{xh\,\mid\, x\in I,\ h\in H\}$
is  $\epsilon$-equidistributed modulo $p$
if $I$ is an interval of $\Fp$ of size $\ge p^{\delta}$
and $H$ is a subgroup of $\Fp^*$
such that $|H|>(\log p)^A$ for some large $A$
depending on the positive real numbers $\epsilon$ and $\delta$. It means that $I\cdot H$ \textit{becomes equidistributed}
when $p$ tends to infinity if one assumes more precisely  that
\begin{equation}\label{pll}
\frac{\log|H|}{\log\log p}\to\infty.
\end{equation}

In Section \ref{s2}, we combine these two Bourgain's statement and obtain a result (cf. Proposition \ref{p1}) on equidistribution for sets $f(I)\cdot H$ where
$H$ is a subgroup of $\Fp^*$,
$I$ is an interval of $\Fp$, $|I|\ge p^{\delta}$ and
$f$ is a non constant polynomial.

Note that $H$ is a coset of a subgroup of $\Fp^{*}$ if and only if  $|H\cdot H|=|H|$.
By relaxing this condition, namely  if the doubling constant of $H$ defined by
$\sigma(H)=|H\cdot H|/|H|$ satisfies $\sigma(H)<K$ ($K>1$), the Green-Ruzsa theorem (i.e. Freiman's theorem in arbitrary abelian groups)  implies that $H$ is well-structured (cf. \cite{GR}).
In this paper, we will focus on the such subsets of $\Fp^{*}$ with doubling constant $\sigma(H)< 2$.  In this restricted case, we will obtain  more elaborated information from the famous Kneser theorem. It suggests the following definition:
we say that a subset of $\Fp^*$ is an \textit{approximate subgroup}
if $|H\cdot H|<2|H|$. It is not difficult to see that
Bourgain-Glibichuk-Konyagin's result (cf. \cite{BGK}) quoted above can be easily extended to approximate subgroup.
In section \ref{s3}, we show our main result
(cf. Theorem \ref{th3}) by extending  Proposition \ref{p1}
to the case
where $H$ is an {approximate subgroup}.

In the last section, we investigate the question of
the existence of residues of a given small subgroup
of $\Fp^*$ in the sumset $A+B$ for two arbitrary
subsets of $\Fp$.

We stress the fact that Bourgain's condition \eqref{pll} on the polylogarithmic size of $H$ is essential in our proofs.
By taking $p=2^q-1$ a Mersenne prime, we can observe
that the multiplicative subgroup $H$ generated by $2$
has cardinality $q =\log(p+1)/\log2> \log p$. Nevertheless, $H$ is not
$\epsilon$-distributed modulo $p$ since $H\cap((p+1)/2,p)=\varnothing$.
Moreover, if $I$ is the interval $(1,2^{\delta q})$ in $\Fp^*$ with $0<\delta<1/2$,
then $I\cdot\{2^j,\ 0\le j\le (1-\delta) q-1\}\subset (0,(p+1)/2)$. This implies that
$|(I\cdot H)\cap((p+1)/2,p)|\le 2^{\delta q}|I|
=o(|I||H|)$, thus $I\cdot H$ is not $\epsilon$-equidistributed when $p$ is large enough
(assuming the Mersenne conjecture
which asserts that there are infinitely many Mersenne
primes, see e.g. \cite{D}).


These questions are related to results and problems
quoted in \cite{C}.

In order to prove our results,
 we will argue by induction on the degree of $f$,
 on the back of Bourgain's result,
using a squaring operation on trignometric sums
and Kneser's theorem on the structure of small \textit{doubling} sets
in Abelian groups.

\section{\bf A result of asymptotic  equidistribution for subgroups of $\Fp^*$}
\label{s2}

Let $H$ be a subset of $\Fp^*$ and $I$ be an interval of $\Fp$, that is
$$
I=\{ax+b\pmod{p}\;\mid\, 0\le x\le |I|-1\},
$$
for some $a\in\Fp^*$ and $b\in\Fp$. We also fix a polynomial $f$
of degree $\ge1$. We consider firstly the question of equidistribution modulo $p$ of the set of residues
$$
f(I)\cdot H:=\{f(z)h\,\mid\, z\in I,\ h\in H\}
$$
 as $p$ tends  to infinity. We prove

\begin{prop}\label{p1}
Let $k$ be positive integer,
$c$ be a positive real number and $\epsilon,\delta\in(0,1]$ be
real numbers. Then there exist $p_0=p_0(k,c,\epsilon,\delta)$,
$A=A(k,c,\epsilon,\delta)$ such that for any prime $p\ge p_0$, any
subgroup $H^*$ of $\Fp^{*}$ with $|H^*|>(\log p)^{A}$, any $u\in\Fp^*$,
any subset $H$ of $uH^*$
with $|H|>c|H^*|$, any interval $I\subset \Fp$ with
$|I|\ge p^{\delta}$ and any $f(x)\in\Fp[x]$ with $\deg f=k$, one has
$$
\Big|\sum_{h\in H}\sum_{z\in I}e_p\big(hf(z)\big)\Big|
\le \epsilon |H||I|.
$$
\end{prop}

For any $a\in\Fp^*$ and any subset $X$ of $\Fp^*$, we denote
\begin{equation}\label{eq0}
N(X,a,\gamma):=\big|\{x\in X\,\mid\, |ax|_p<p^{1-\gamma}\}\big|,
\end{equation}
where $|x|_p$ means the unique nonnegative integer less than $p/2$ congruent
to $|x|$ modulo $p$.

Let $\gamma>0$.
Our aim is  to use Bourgain's result on the distribution of $g^n$ modulo $p$, $n\ge0$, where $g\in\Fp^*$ is fixed
(cf. \cite{B}).    It implies that if  the size of $H^*$  is  sufficiently  large,  namely
$$
\frac{\log|H^*|}{\log\log p} > A_1,
$$
where $A_1$ is a computable large constant in terms of $\epsilon$, $c$ and $\gamma$,  then one has
$$
\max_{a\in\Fp^*}N(H^*,a,\gamma)\le \frac{\epsilon c|H^*|}2,
$$
for any sufficiently large prime number $p$.

By assumption on $H$, we  deduce that  $ \max_{a\in\Fp^*}N(H,a,\gamma)\le \frac{\epsilon |H|}2$.

We assume that $|I|\ge p^{\delta}$. Then for any integer $r\in(1,p-1)$, we have
$$
T_r:= \Big|\sum_{h\in H}\sum_{z\in I}e_p\big(rf(z)h\big)\Big|\le\epsilon|I||H|.
$$
Indeed, we get for $f(z)=az+v$, $a\in\Fp^*$ and $v\in \Fp$,  
\begin{equation}\label{eq1}
T_r\le
\sum_{h\in H}
\Big|\sum_{z\in I}e_p\big(razh\big)\Big|
\le
 \sum_{h\in H}\min(|I|,\|arh/p\|^{-1})
\le |I|N(H,ar,\gamma)+p^{\gamma}|H|\le \epsilon|I||H|,
\end{equation}
if one chooses $\gamma =\delta/2$ and if $p$ is large enough. By letting $A(1,c,\epsilon,\delta)=A_1$ we get
the result for $k=1$.

For a general non constant polynomial $f(x)\in\mathbb{Z}[x]$, we argue by induction on $k\ge1$. For $k=1$,
it has been done above.
Assume now that the property holds for some $k\ge1$. Let $f$ be
of degree $k+1$. By letting
$$
S(h)=\sum_{z\in I}e_p(hf(z)),
$$
we have
$$
\sum_{h\in H}|S(h)|^2
=\sum_{h\in H}\sum_{x,y\in I}e_p(h(f(x)-f(y)))
\le|H||I|+2\sum_{z\in I}\Big|\sum_{h\in H}\sum_{y\in I}
e_p(h g_z(y))\Big|,
$$
where $g_z(y):=f(y+z)-f(y)$ is a polynomial of degree $k$.
By the inductive hypothesis, we get
$$
\sum_{h\in H}|S(h)|^2\le 3\epsilon |H||I|^2.
$$
Thus the set
$$
H':=\{h\in H\,|\, |S(h)|>\epsilon^{1/3}|I|\}
$$
has cardinality satisfying
$$
\epsilon^{2/3}|H'||I|^2<3\epsilon|H||I|^2,
$$
yielding $|H'|< 3\epsilon^{1/3}|H|$. It follows that
\begin{align*}
\Big|\sum_{h\in H}S(h)\Big|&\le
\sum_{h\in H'}|S(h)|+\sum_{h\in H\smallsetminus H'}|S(h)|\\
&\le |H'||I|+|H|\epsilon^{1/3}|I| <  4\epsilon^{1/3}|I||H|.
\end{align*}
The result is proved.

We can derive from the proof that Proposition \ref{p1}
holds uniformly
for any  polynomial of degree less than $k(\epsilon)=\frac{\log\log(1/\epsilon)}{\log3}$.

\section{\bf Extension to approximate multiplicative subgroups}\label{s3}

We recall that an \textit{approximate} subgroup of $\Fp^{*}$ is any subset
$H$ of $\Fp^{*}$ such that $|H\cdot H|<2|H|$.
By Kneser's Theorem, we get the following structure for such
a subset:

\begin{lemma} Let $\eta>0$ and $H\subset \Fp^*$. If $|H\cdot H|<(2-\eta)|H|$,
then there exist a positive  integer $m\le1/\eta$, a subgroup $H^*$ of $\Fp^{*}$ and $u_1,\dots,u_m\in H$ such that
$$
H\subset \bigcup_{i=1}^m u_iH^{*}.
$$
and
$$
\frac{|H|}{m}\le|H^*|\le  \frac{2-\eta}{2m-1}|H|.
$$
\end{lemma}

We can now generalize Bourgain's result to
approximate subgroup multiplying by  the image of an interval
by a polynomial.

\begin{thm}\label{th3}
Let $H$ be an approximate subgroup of $\Fp^{*}$ with size larger than polylogarithmic in $p$ and
$f$ be a polynomial. Then for any interval $I$ in $\Fp$ of size $p^{\delta}$,
$f(I)\cdot H$ is equidistributed modulo $p$ as $p$ tends to infinity.
\end{thm}


Let $\epsilon,\delta,\eta>0$ be positive real numbers and $k$ be a positive  integer. We assume that $|H\cdot H|\le (2-\eta)|H|$ and $|H|>(\log p)^{B}/\eta$ where
$B=A(k,\epsilon\eta,\epsilon,\delta)$ is defined in Proposition \ref{p1}. By the previous lemma, we may write
$$
H=\bigcup_{i=1}^mH_i
$$
where $H_i=u_iH^{*}\cap H$, $i=1,\dots,m$. Let
$$
\Lambda=\{i\le m\,\mid\,
|H_i|\le \epsilon |H^*|/m\},
$$
and $\Lambda'=\{1,2,\dots,m\}\smallsetminus \Lambda$.

For $j\in \Lambda'$, we have both $|H_j|>\epsilon|H^{*}|/m
\ge \epsilon\eta|H^*|$
and $|H^*|\ge |H|/m>(\log p)^{B}$.
Since Proposition \ref{p1} holds for cosets of a multiplicative subgroup as well, we obtain
\begin{align*}
\Big|\sum_{h\in H}\sum_{z\in I}e_p\big(hf(z)\big)\Big|&\le
\sum_{i=1}^m\Big|\sum_{h\in H_i}\sum_{z\in I}e_p\big(hf(z)\big)\Big|\\
&\le
\sum_{i\in \Lambda}|H_i||I|+\sum_{i\in \Lambda'}\epsilon |I||H_i|\\
&\le \epsilon|\Lambda||I||H^{*}|+\epsilon|I||H|\le 3\epsilon|I||H|.
\end{align*}

\section{\bf Remarks}

It is worth mentioning that a close question related to multiplicative subgroups of $\Fp^*$ can be considered: does the equation $a+b=h$, $(a,b,h)\in A\times B\times H$ be solvable for
any subsets $A,B$ of $\Fp$ and any subgroup of $\Fp^*$ ?
Of course, $A,B$ and $H$ must be large enough in terms of
$p$. This type of question takes its origin in \cite{Sa}
and has been hugely investigated since (see e.g.
\cite{S}, \cite{Sh} and \cite{HH}).

By the use of Fourier analysis in $\Fp^*$ with
ingredients of \cite{S} (see also \cite{Sh}), it can be shown that it is the case
if
\begin{equation}\label{eqn4}
|A|>p^{\epsilon},\quad |B|>p^{1/2+\epsilon},
\quad |H|>p^{1-\delta},
\end{equation}
where $\delta=\delta(\epsilon)>0$. The proof runs as follows.
The number of solutions of the equation $a+b=h$ is equal
to
$$
N=\frac{|A||B||H|}{p-1}+\frac{1}{p-1}\sum_{r=1}^{p-2}
\sum_{(a,b)\in A\times B}\chi_r(a+b)\sum_{h\in H}\overline{\chi_r}(h),
$$
where $\chi_r$ denotes the multiplicative character modulo $p$ defined by
$$
\chi_r(x)=\exp\left(\frac{2\pi i r\  \text{ord}(x)}{p-1}\right),
$$
and $\text{ord}(x)$ denotes the discrete logarithm of $x$ in base
$g$ for some fixed primitive root $g$ modulo $p$.
The summation on $h$ is $|H|$ or $0$ according to the
fact that $r$ divides $|H|$ or not. Hence
$$
N=\frac{|A||B||H|}{p-1}+\frac{|H|}{p-1}\sum_{s=1}^{(p-1)/|H|-1}
\sum_{(a,b)\in A\times B}\chi_{s|H|}(a+b).
$$
By Shparlinski's result (cf. eq. 14 in \cite{S}), the summation on
$(a,b)$ is $O(|A||B|p^{-\delta'})$ for any $s$, hence $N>0$
by \eqref{eqn4} if we consider $\delta<\delta'$ and $p$ sufficiently large.

The same result with a stronger assumption on $|A|$ and
$|B|$ and by relaxing the one on $|H|$ is a consequence
of the corollary  to Theorem 2.4 of \cite{H}:
$$
\text{$a+b=h$ is solvable if}\quad |A||B|>p^{2-\epsilon},\quad |H|>p^{1/3+\delta},
$$
where $\epsilon\to0$ as $\delta\to0$.


\begin{thebibliography}{9}


\bibitem{B0} Bourgain, J.; Mordell's exponential sum estimate revisited.  \textit{J. Amer. Math. Soc.}  \textbf{18}  (2005),  no. 2, 477--499.


\bibitem{B} Bourgain, J.;
On the distribution of the residues of small multiplicative
subgroups of $\Fp$. \textit{Israel J. of Math.} \textbf{172} (2009), 61--74.

\bibitem{BGK} Bourgain, J.; Glibichuk, A.; Konyagin S.;
Estimate for the number of sums and products and for exponential sums in fields of prime order.
\textit{J. London Math. Soc.}
\textbf{73} (2006), 380--398.

\bibitem{C} Chang M.-C.; Some problems in combinatorial number theory.  \textit{Integers}  \textbf{8}  (2008),  no. 2, A1, 11 pp.

\bibitem{D} Dickson L.E.; History of the theory of numbers. Chelsea Publishing, New York, 1971.

\bibitem{GR}Green B.J; Ruzsa I.Z.;  Freiman's theorem in an arbitrary abelian group.
\textit{J. London Math. Soc.} \textbf{75} (2007), 163--175.

\bibitem{H} Hegyv\'ari N.; Some remarks on multilinear sums
and their applications. Preprint, 2010.

\bibitem{HH} Hegyv\'ari N.; Hennecart F.; Explicit construction of extractors and expanders,
{\it Acta Arith.} {\bf140} (2009), 233--249.

\bibitem{Sa} S\'ark\"ozy, A.; On sums and products of residues
modulo $p$. \textit{Acta Arith.} {\bf118} (2005), 403--409.

\bibitem{Sh} Shkredov I.D., On monochromatic solutions of some non linear equations, preprint (2009).

\bibitem{S} Shparlinski I.E.; On the solvability of bilinear equations in finite fields.
\textit{Glasg. Math. J.}  \textbf{50}  (2008),  no. 3, 523--529.

\bibitem{TV} Tao T; Vu V.H.; Additive combinatorics. Cambridge Studies in Advanced Mathematics, 105. Cambridge University Press, Cambridge, 2006. xviii+512 pp.





\end{thebibliography}
\end{document}